\documentclass[10pt,reqno]{amsart}

\usepackage{amsmath,amsfonts,amsthm,amssymb,amsxtra,dsfont}

\newtheorem{theorem}{Theorem}
\newtheorem{proposition}{Proposition}

\theoremstyle{definition}

\numberwithin{equation}{section}

\newcommand\1{{\ensuremath {\mathds 1} }}

\newcommand{\C}{\mathbb{C}}

\renewcommand{\epsilon}{\varepsilon}

\newcommand{\R}{\mathbb{R}}

\renewcommand{\leq}{\leqslant}
\renewcommand{\geq}{\geqslant}

\renewcommand{\to}{\rightarrow}

\begin{document}

\title[Ground States for the boson star equation]{Uniqueness of ground states \\ for the $L^2$-critical boson star equation}

\author{Rupert L. Frank}
\address{Rupert L. Frank, Dept.~of Mathematics, Princeton University, Washington Road, Princeton, NJ 08544, USA}
\email{rlfrank@math.princeton.edu}

\author{Enno Lenzmann} 
\address{Enno Lenzmann, Dept.~of Mathematics, MIT, 77 Massachusetts Avenue, Cambridge, MA 02139, USA}
\email{lenzmann@math.mit.edu}


\begin{abstract}
We establish uniqueness of ground states $u(x) \geq 0$ for the $L^2$-critical boson star equation 
$$\sqrt{-\Delta} \, u - ( |x|^{-1} \ast |u|^2 ) u = -u, \quad \mbox{in $\R^3$.}$$ 
The proof blends variational arguments with the harmonic extension to the halfspace $\R^4_+$. Apart from uniqueness, we also show radiality of ground states (up to translations) and the nondegeneracy of the linearization. Our results provide an indispensable basis for the blowup analysis of the time-dependent $L^2$-critical boson star equation. The uniqueness proof can be generalized to different fractional Laplacians $(-\Delta)^s$ and space dimensions.
\end{abstract}

\maketitle

\section{Introduction and Main Results}

Fractional powers of the Laplacian arise in a numerous variety of equations in mathematical physics; see, e.\,g., \cite{CaSi,AbBoFeSa,MaMcTa,LiYa,ElSc,FrLe} and references therein. Here, a fundamental r\^ole is often played by  {\em ground states,} by which we mean nonnegative solutions $u(x) \geq 0$ in $H^s(\R^d)$ that satisfy an equation of the form
\begin{equation}
\label{eq:general}
(-\Delta)^s u + F(u) = 0 , \quad \mbox{for $u : \R^d \to [0,\infty)$}.
\end{equation}
As usual, the fractional Laplacian $(-\Delta)^s$, with $0<s<1$, is defined by its multiplier $|\xi|^{2s}$ in Fourier space, and $F(u)$ denotes some given (not necessarily local) nonlinearity. 

In contrast to the question of existence, it is fair to say that extremely little is known about the uniqueness (modulo symmetries) of ground states for problems like \eqref{eq:general}, since classical shooting and ODE methods are clearly not applicable. To the authors' knowledge, there are only two cases where uniqueness of ground states for fractional equations is known to hold:
\begin{enumerate}
\item[(i)] Solitary waves for the Benjamin-Ono and ILW equations in $d=1$ dimension; see \cite{AmTo,AlTo}.
\item[(ii)] Optimizers for the Hardy-Littlewood-Sobolev (HLS) inequality in dimensions $d \geq 1$; see \cite{ChLiOu,Li}.
\end{enumerate}
However, both results rest on a very specific feature of each problem: In case (i), the proof is intimately linked to {\em complex analysis} and subtle identities that can only be applied to the Benjamin-Ono and ILW equations; whereas the case (ii) exhibits the {\em conformal symmetry} of the HLS inequality, which plays a key role in the uniqueness proof.


In this note, we establish uniqueness and radiality of ground states for the {\em $L^2$-critical boson star equation} 
\begin{equation}
\label{eq:eq}
 \sqrt{-\Delta} \, u - \big ( |x|^{-1} \ast |u|^2 \big ) u   =  -u \quad \mbox{in $\R^3$},
\end{equation}
where $u \in H^{1/2}(\R^3)$, $u \geq 0$ and $u \not \equiv 0$; and $\ast$ denotes convolution. The uniqueness proof will be sketched below and detailed in \cite{FraLe}. The main idea employs variational arguments for a certain nonlinear boundary-value problem in the halfspace $\R^4_+ = \R^3 \times \R_+$. Indeed, we do think that some key ideas will be beneficial for the uniqueness proof of ground states for problems like \eqref{eq:general} in more generality; see also \cite{FraLe} for extensions.

Equation \eqref{eq:eq} plays a central role in the mathematical theory of {\em boson stars,} initiated in \cite{LiYa}. Recently, it was shown in \cite{ElSc} that the {\em effective dynamics} of boson stars can be described in terms of the evolution equation (with mass parameter $m \geq 0$)
\begin{equation}
\label{eq:boson}
i \partial_t \psi = \sqrt{-\Delta + m^2} \, \psi - \big ( |x|^{-1} \ast |\psi|^2 \big ) \psi
\end{equation}
where $\psi : [0,T) \times \R^3 \to \C$ is a wave field. In fact, this dispersive nonlinear $L^2$-critical PDE displays stable/unstable traveling solitary waves and finite-time blowup, reflecting various physical situations; see \cite{FrLe,FrJoLe,FrJoLe2,Le2}. In particular, the ground states for \eqref{eq:eq} are of fundamental interest, since the (unstable) solitary waves $\psi(t,x) = e^{it} u(x)$ provide {\em nondispersive, critical objects} for equation \eqref{eq:boson} in the $L^2$-critical regime when $m=0$. Furthermore, as shown in \cite{Le2} for all $m\geq 0$, the solutions $\psi(t,x)$ in $H^{1/2}(\R^3)$ with $L^2$-mass $\| \psi(t) \|_{2}^2 < \| u \|_2^2$ extend globally in time; whereas solutions with $\| \psi(t) \|_{2}^2 > \| u \|_2^2$ can lead to {\em finite-time blowup,} as proven in \cite{FrLe} for radial solutions. Thus our main results here form an indispensable basis for the analysis and classification of blowup for  equation \eqref{eq:boson} for initial data close to $e^{it} u(x)$.

\medskip
\begin{theorem}\label{th:main} 
Equation \eqref{eq:eq} has a unique ground state up to translations. More precisely, there exists a unique, radial and positive function $Q =Q (|x|) \in H^{1/2}(\R^3)$ such that, for any ground state $u(x) \geq 0$ solving \eqref{eq:eq}, there exists $a\in\R^3$ such that $u(x)=Q(|x-a|)$. 

Moreover, the function $Q(|x|)$ is strictly decreasing and analytic; and its Fourier transform $\hat{Q}(|\xi|)$ is positive, radial and non-increasing.
\end{theorem}



\medskip
In \cite{FraLe} we generalize Theorem \ref{th:main} to dimensions $d \geq 3$, convolution kernels $|x|^{2-d}$, and fractional Laplacians $(-\Delta)^s$. We remark that the proof of radiality is based on the methods of {\em moving planes} and a non-local Hopf-type lemma; see \cite{FraLe}. See also the moving plane method developed in \cite{ChLiOu} and extended in \cite{MaZh} as an alternative approach to prove radiality. The analyticity of $Q$ follows from  $\| |\xi|^n \hat Q \|_\infty \leq (\mathrm{const} \cdot n)^n $, which follows by extending the method of \cite{LiBo} to systems of non-local elliptic equations.

\medskip
Our second main result concerns the so-called {\em nondegeneracy} of the linearization of \eqref{eq:eq} around the unique radial ground state $Q = Q(|x|) \in H^{1/2}(\R^3)$. More precisely, this means that the linear operator
\begin{equation}
\label{eq:Lplus}
L_+ \xi = \sqrt{-\Delta} \,  \xi + \xi - \big (|x|^{-1} \ast |Q|^2 \big) \xi - 2 Q \big ( |x|^{-1} \ast (Q\xi) \big )
\end{equation}
has a kernel that is entirely due to the translational invariance of equation \eqref{eq:eq}. 

\medskip
\begin{theorem} \label{th:kernel}
Let $Q = Q(|x|) \in H^{1/2}(\R^3)$ denote the unique radial solution of \eqref{eq:eq} given by Theorem \ref{th:main} above. Define the operator $L_+$ as in \eqref{eq:Lplus} acting on $L^2(\R^3)$ with domain $H^1(\R^3)$. Then $L_+$ is nondegenerate, i.\,e., its kernel satisfies 
$$
\ker L_+ = \mathrm{span} \, \{ \partial_{x_1} Q, \partial_{x_2} Q, \partial_{x_3} Q \} \,.
$$
\end{theorem}

\medskip
In the context of local nonlinear Schr\"odinger equations (NLS), the nondegeneracy of ground states is a well-known fact; see \cite{ChGuNaTs,We}. However for (NLS), the proof involves some Sturm-Liouville theory, which clearly is not applicable to $L_+$ given by \eqref{eq:Lplus} due to its nonlocality. To overcome these difficulties, we adapt heat kernel arguments (developed in \cite{Le}) and use insights from the proof of Theorem \ref{th:main} given in \cite{FraLe}.

\section{Sketch of Uniqueness Proof for Radial Solutions}

We remind the reader that complete proofs will be given in \cite{FraLe}.

For any radial ground state $u =  u(|x|) > 0$ of \eqref{eq:eq}, we consider its harmonic extension $U=U(x,t)$ onto the halfspace $\R^4_+ = \{ (x,t) \in \R^3 \times \R_+ \}$ given by $U(x,t):=(\mathcal{P}_t \ast u)(x)$, where $\mathcal{P}_t(z)$ denotes the Poisson kernel in three dimensions. Since $u \in H^{1/2}(\R^3)$ and by the spatial decay $u(x) \lesssim \langle x \rangle^{-4}$, we deduce that $U \in H^1(\R^4_+)$ holds. In summary, we see that $U \in H^1(\R^4_+)$ solves the nonlinear  boundary-value problem:
\begin{equation} \label{eq:halfspace} -\Delta U = 0  \quad \mbox{in $\R^4_+$}, \quad \mbox{and} \quad - \partial_t U + (\Phi_u - 1) U = 0 \quad \mbox{on $\partial \R^4_+ = \R^3 \times \{ 0 \}$}. 
\end{equation}
Here we applied (for later purpose) Newton's theorem to express the convolution in \eqref{eq:eq}, where $\Phi_u(x) := \int_{|y| < |x|} (|y|^{-1} - |x|^{-1}) |u(y)|^2 \, dy$. The use of Newton's theorem is inspired by Lieb's uniqueness proof in \cite{Lie}. Also, we rescaled $u \mapsto e^{-3/2} u(e^{-1} \cdot )$ for some $e > 0$ in \eqref{eq:halfspace}. Next, we introduce the quadratic form
\begin{equation}
\mathcal{A}_u[\psi] := \int \! \! \int_{\R^4_+} |\nabla_{(x,t)} \psi|^2 \,dx \,dt + \int_{\R^3} (\Phi_u - 1) |\psi(x,0)|^2 \, dx ,
\end{equation}
for any $\psi \in H^1(\R^4_+)$. As an important fact, we can derive the following lower bound.

\medskip
\begin{proposition} \label{prop:piepenbrink} 
We have $\mathcal{A}_u[\psi] \geq 0$ for all $\psi \in H^1(\R^4_+)$. Moreover, we have $A_u[\psi] = 0$ if and only if $\psi \equiv \lambda U$ for some constant $\lambda \in \mathbb{C}$.
\end{proposition}
\medskip
Since $U(x,t) > 0$  and $A_u[U] = 0$ by inspection, we remark that such a lower bound for $A_u[\psi]$ is reminiscent of {\em Allegretto-Piepenbrink theory:} Positive supersolutions imply nonnegativity of Schr\"odinger operators. We now proceed to the heart of the uniqueness proof. Suppose that $u = u(|x|) >0$ and $v = v(|x|) > 0$ are two radial and positive solutions to \eqref{eq:eq} such that $u \not \equiv v$. 

Using \eqref{eq:eq} and the positivity of $u$ and $v$, it is easy to see that the two solutions necessarily intersect. Further, by analyticity arguments, there is a smallest $R > 0$ such that $u(|x|) = v(|x|)$ when $|x| = R$. We can assume $u(|x|) > v(|x|)$ for $0 < |x| < R$. Next, let $U$ and $V$ be the harmonic extensions of $u$ and $v$, respectively, and let $\mathcal{A}_u[\psi]$ and $\mathcal{A}_v[\psi]$ be the corresponding quadratic forms. The idea is now to violate the lower bound given by Proposition \ref{prop:piepenbrink}. Indeed, as detailed in \cite{FraLe}, we find a set $\Omega \subset \overline{\R^4_+}$ such that $W := (U-V) \1_\Omega $ satisfies: $W \in H^1(\R^4_+)$, $W \geq 0$, $W \not \equiv 0$, and $-\Delta W = 0$ in $\R^4_+ \setminus \partial \Omega$, as well as
$$
-\partial_t W + \frac{1}{2}(\Phi_u + \Phi_v) W = W - f \quad \mbox{on} \quad  \{ |x| < R \} \times \{ t=0 \}. 
$$
Moreover, we have $W \equiv 0$ on $\{ |x| \geq R \} \times \{ t =0 \}$ and $f = \frac{1}{2} (\Phi_u - \Phi_v) (u+v)$. (Note that $W$ is {\em not} the harmonic extension of $(u-v)\1_{|x| < R}$.) Integrating by parts and using that $u > v > 0$ and $\Phi_u > \Phi_v$ for $0 < |x| < R$, we find that $\mathcal{A}_u[W] + \mathcal{A}_v[W] = -2 \int_{|x| < R} f W |_{t=0} < 0$, contradicting Proposition \ref{prop:piepenbrink}.

\section*{Acknowledgments}
R.L.~acknowledges support through DFG grant FR 2664/1-1, NSF grant PHY 06 52854; and E.L.~is partially supported by  NSF grant DMS-070249 and CTS at ETH Z\"urich. We thank Joachim Krieger, Mathieu Lewin and Pierre Rapha\"el for valuable discussions.

\end{document}